\documentclass{article}

\usepackage[utf8]{inputenc}
\RequirePackage[OT1]{fontenc}

\usepackage{mathptmx} 

\usepackage[english]{babel}

\usepackage{amsmath,amscd,amsfonts,amsthm}
\usepackage{mathrsfs}
\usepackage[psamsfonts]{amssymb}
\usepackage{mathtools}
\usepackage{dsfont}
\RequirePackage[numbers]{natbib}
\RequirePackage[colorlinks,citecolor=blue,urlcolor=blue]{hyperref}

\numberwithin{equation}{section}
\theoremstyle{plain}
\newtheorem{thm}{Theorem}[section]

\newtheorem{lem}{Lemma}
\newtheorem{prop}{Proposition}

\newcommand{\E}{\mathbb{E}}
\newcommand{\tr}{\mathop{\mathrm{Tr}}}
\newcommand{\ft}{\hat{f}_t}
\newcommand{\pen}{\mathop{\mathrm{pen}}}
\newcommand{\KL}{\mathrm{KL}}
\newcommand{\price}{\mathop{\mathrm{price}}}

\title{PAC-Bayesian aggregation of affine estimators}

\author{Lucie Montuelle\thanks{RTE, La D\'efense, France,
    \url{lucie.montuelle@rte-france.com}} and
    Erwan Le Pennec\thanks{CMAP/XPOP, \'Ecole Polytechnique, France, \url{erwan.le-pennec@polytechnique.edu}}}

  \date{January 2018}

\begin{document}

\maketitle

\begin{abstract}
 Aggregating estimators using exponential weights
depending on their risk appears optimal in expectation
but not in probability. We use here a slight
overpenalization to obtain oracle inequality in probability for 
such an explicit aggregation procedure. We focus on the fixed design
regression framework and the aggregation of affine estimators and
obtain results for a large family of affine estimators under a
non necessarily independent sub-Gaussian noise assumptions.
\end{abstract}

\tableofcontents

\section{Introduction}
\label{sec:intro}

We consider here a classical fixed design regression model \[\forall
i\in \{1,\ldots, n\},\,  Y_i = f_0(x_i) + W_i \] with $f_0$ an unknown
function, $x_i$ the fixed design points and $W=(W_i)_{i\leq n}$ a
centered sub-Gaussian noise. We assume that we have at hand a family
of affine estimate $\{\ft(Y)=A_t Y + b_t |A_t \in \mathcal{S}^+_n(\mathbb{R}), b_t \in \mathbb{R}^n, t \in
\mathcal{T}\}$, for instance a family of projection estimator, of
linear ordered smoother in a basis or in a family of basis. The most
classical way to use such a family is to select one of the estimate
according to the observations, for instance using a penalized empirical risk
principle. A better way is to combine linearly those estimates with
weights depending of the observation. A simple strategy
is  the Exponential
Weighting Average in which all those estimate are averaged with a
weight proportional to $\exp\left( -\frac{\widetilde{r_t}}{\beta}\right) \pi(t)$ where
$\widetilde{r_t}$ is a (penalized) estimate of the risk of
$\ft$. This strategy is not new  nor optimal as explained
below but is widely used in practice. In this article, we analyze the
performance of this simple EWA estimator by providing oracle
inequalities in probability under mild sub-Gaussian assumption on the
noise.

Our aim is to obtain the  best possible estimate of the function $f_0$
at the grid points. This setting is probably one of the most common in
statistics and many regression estimators are available in the
literature. For non parametric estimation, Nadaraya-Watson estimator
\citep{Nadaraya1965,MR0185765} and its fixed design counterpart
\citep{MR767241} are widely used, just like projection estimators
using trigonometric, wavelet \citep{MR1323344} or spline
\citep{MR1045442} basis for example. In the parametric framework,
least squares or maximum likelihood estimators are commonly employed,
sometimes with minimization constraints, leading to
\textsc{lasso}~\citep{Tibshirani94regressionshrinkage},
ridge~\citep{MR2722294}, elastic net~\citep{MR2137327},
\textsc{aic}~\citep{MR0483125} or \textsc{bic} \citep{MR0468014}
estimates.

Facing this variety, the statistician may wonder which procedure
provides the best estimation. Unfortunately, the answer depends on the
data. For instance, a rectangular function is  well approximated by
wavelets but not by trigonometric functions.
Since the best estimator is not known in advance, our aim is to mimic
its performances in term of risk. This is theoretically guaranteed by
an oracle inequality:
\[R(f_0,\tilde{f}) \leq C_n \inf_{t\in \mathcal{T}} R(f_0,\ft)+
\epsilon_n\]
comparing the risk of the constructed estimator $\tilde{f}$ to the
risk of the best available procedure in the collection
$\{\ft, t \in \mathcal{T}\}$.
Our strategy is based on convex combination of these preliminary
estimators and relies on PAC-Bayesian aggregation to obtain a single
adaptive estimator. We focus on a wide family, commonly used in
practice : affine estimators
$\{\ft(Y)=A_t (Y-b) + b + b_t |A_t \in \mathcal{S}^+_n(\mathbb{R}), b_t \in \mathbb{R}^n, t \in
\mathcal{T}\}$ with $b \in \mathbb{R}^n$ a common recentring.

Aggregation procedures have been introduced by \citet{Vovk1990,
  Littlestone1994212, Cesa-Bianchi1997, MR1765620}. They are a central
ingredient of bagging~\citep{raey}, boosting
\citep{MR1348530,Schapire1990} or random forest (\citet{Amit1997} or
\citet{randomforests}; or more recently \citet{MR2447310, MR2719877,
  MR2930634, Genuer2011}).

The general aggregation framework is detailed in \citet{MR1775640} and
studied in \citet{MR2163920,MR2483528} through a PAC-Bayesian framework
as well as in \citet{MR1762904, MR1790617, MR1804557, MR1946426,
  MR1997174, MR2028357, MR2044592}. See for instance \citet{MR2503001}
for a survey. Optimal rates of aggregation in regression and density
estimation are studied by \citet{Tsybakov2003, MR2356820, MR2356821,
  Rigollet2006} and \citet{MR2364224}.

A way to translate the confidence of each preliminary estimate is to
aggregate according to a measure exponentially decreasing when the
estimate's risk rises. This widely used strategy is called
exponentially weighted aggregation. More precisely, as explained
before, the weight of each
element $\ft$ in the collection is proportional to
$\exp\left( -\frac{\widetilde{r_t}}{\beta}\right) \pi(t)$ where
$\widetilde{r_t}$ is a (penalized) estimate of the risk of
$\ft$, $\beta$ is a positive parameter, called the temperature, that
has to be calibrated and $\pi$ is a prior measure over
$\mathcal{T}$. The main interest of exponential weights resides in
Lemma~\ref{Gibbs} \citep{MR2483528} since they explicitly minimize the
aggregated risk penalized by the Kullback-Leibler divergence to the prior measure $\pi$.  Our aim is to give sufficient conditions on the
risk estimate $\widetilde{r_t}$ and the temperature $\beta$ to obtain
an oracle inequality for the risk of the aggregate. 
Note that when the family $\mathcal{T}$ is countable, the exponentially
weighted aggregate is a weighted sum of the preliminary estimates.

This procedure has shown its efficiency, offering lower risk than
model selection because we bet on several estimators. Aggregation of
projections has already been addressed by~\citet{MR2242356}. They have
proved by the mean of an oracle inequality, that in expectation, the
aggregate performs almost as well as the best projection in the
collection.  Those results have been extended to several settings and
noise conditions \citep{Dalalyan2007, DT08, MR2926142, MR2543587,
  Dalalyan_Hebiri_Meziani_Salmon13, MR2860324, MR2981354, MR3025131,
  MR2999166, MR3025134} under a \emph{frozen} estimator assumption:
they should not depend on the observed sample. This restriction, not
present in the work by \citet{MR2242356}, has been removed by
\citet{MR3059085} within the context of affine estimator and
exponentially weighted aggregation. 
Nevertheless, they make additional assumptions on the matrices $A_t$
and the Gaussian noise to obtain an optimal oracle inequality in
expectation for affine estimates. Very sharp results have been obtained in 
\citet{golubev12:_expon}, \citet{chernousova13:_order} and
\citet{golubev14:_concen_inequal_expon_weigh_method}. Those papers,
except the last one, study a risk in expectation.

Indeed, the Exponential Weighting Aggregation is not optimal anymore
in probability.  \citet{MR3015047} have indeed proved 
the sub-optimality in deviation
of exponential weighting, not allowing to obtain a sharp oracle
inequality in probability. Under strong assumptions and independent
noise, \citet{bellec} provides a sharp oracle inequality with optimal
rate for another aggregation procedure called Q-aggregation. It is
similar to exponential weights but the criterion to minimize is
modified and the weights no longer are explicit. Results for the
original EWA scheme exists nevertheless but with a constant strictly
larger than $1$ in the oracle inequality. \citep{MR3192554} obtain
for instance a result under a Gaussian white noise assumption by 
penalizing the risk in the weights and taking a temperature at least
20 times greater than the noise
variance. \citet{golubev14:_concen_inequal_expon_weigh_method} does
not use an overpenalization but assume some ordered structure on the
estimate to obtain a result valid even for low temperature.
 An unpublished work,
by~\cite{gerchinovitztel00653550}, provides also weak oracle inequality
with high probability for projection estimates on non linear
models. \citet{MR2786484} consider \emph{frozen} and bounded
preliminary estimators and obtain a sharp oracle inequality in
deviation for the excess risk under a sparsity assumption, if the
regression function is bounded, with again a
 modified version of exponential weights.

 In this article, we will play on both the temperature and the
 penalization. We will be able to obtain oracle inequalities for the
 Exponential Weighting Aggregation under a general sub-Gaussian noise
 assumption that does not require a coordinate independent setting. We
 conduct an analysis of the relationship between the choice of the
 penalty and the minimal temperature. In particular, 
we show that there is a continuum between the usual noise based
penalty and a sup norm type one
allowing a \emph{sharp} oracle inequality.

\section{Framework and estimate}
\label{sec:Framework}

Recall that we observe
\[ \forall i \in \{1,\ldots, n\},\, Y_i = f_0(x_i) + W_i \]
with $f_0$ an unknown function and $x_i$ the fixed grid points. Our
only assumption will be on the noise. We do not assume any independence
between the coordinates $W_i$ but only that 
$W=(W_i)_{i\leq n} \in \mathbb{R}^n$ is a centered sub-Gaussian
variable. More precisely, we assume that
$\E(W)=0$ and there exists $\sigma^2 \in \mathbb{R}^+$
such that
\[ \forall \alpha \in \mathbb{R}^n,\, \E \left[ \exp \left(\alpha^\top
    W \right) \right] \leq \exp \left(\frac{\sigma^2}{2}
  \|\alpha\|^2_2 \right), \]
where $\|.\|_2$ is the usual euclidean norm in $\mathbb{R}^n$.  If $W$
is a centered Gaussian vector with covariance matrix $\Sigma$ then
$\sigma^2$ is nothing but the largest eigenvalue of $\Sigma$.

The quality of our estimate will be measured through its error at the
design points. More precisely, we will consider the classical
euclidean loss, related to the squared norm
\[ \|g\|^2_2= \sum_{i=1}^n g(x_i)^2.\]
Thus, our unknown is the vector $(f_0(x_i))_{i=1}^n$ rather than the
function $f_0$.

As announced, we will consider affine estimators $\ft(Y)= A_t (Y-b) +
b + b_t$
corresponding to affine smoothed projection.

We will assume that
\[ \ft(Y) = A_t (Y-b) + b + b_t =\sum_{i=1}^n \rho_{t,i} \langle Y-b, g_{t,i} \rangle
g_{t,i} + b + b_t \]
where $(g_{t,i})_{i=1}^n$ is an orthonormal basis,
$(\rho_{t,i})_{i=1}^n$ a sequence of non-negative real numbers and
$b_t \in \mathbb{R}^n$. By construction, $A_t$ is thus a symmetric
positive semi-definite real matrix. 
 We assume furthermore that the matrix collection
 $\{ A_t \}_{t \in \mathcal{T}}$ is such that
$\sup_{t\in \mathcal{T}} \|A_t\|_2 \leq 1$.
For sake of simplicity, we only use the
notation $\ft(Y)= A_t (Y-b) + b + b_t$ in the following. 

To define our estimate from the collection
$\{\ft(Y)=A_t Y+b_t|A_t \in \mathcal{S}^+_n(\mathbb{R}), b_t \in \mathbb{R}^n, t \in
\mathcal{T}\}$,
we specify the estimate $\widetilde{r_t}$ of the (penalized) risk of the estimator
$\ft(Y)$, choose a prior probability measure $\pi$ over $\mathcal{T}$
and a temperature $\beta>0$. We define the exponentially weighted
measure $\rho_{EWA}$, a probability measure over $\mathcal{T}$, by
\begin{align*}
  d\rho_{EWA}(t) = \frac{\exp\left(-\frac{1}{\beta}\widetilde{r_t} \right)}
{\int \exp\left(-\frac{1}{\beta}\widetilde{r_{t'}} \right)
  d\pi(t')} d\pi(t)
\end{align*}
and the exponentially weighted aggregate $f_{EWA}$ by
$f_{EWA}=\int \ft \, d\rho_{EWA}(t)$.
If $\mathcal{T}$ is countable then
\begin{align*}
  f_{EWA} = \sum_{t \in \mathcal{T}}
 \frac{ e^{-\widetilde{r_t}/\beta} \pi_t}{\sum_{t' \in \mathcal{T}} e^{-\widetilde{r_{t'}}/\beta} \pi_{t'}} \ft.
\end{align*}

This construction naturally favors low risk estimates.  
When the temperature goes to zero, this estimator becomes very similar
to the one minimizing the risk estimate while it becomes an
indiscriminate average when $\beta$ grows to infinity. The choice of
the temperature appears thus to be crucial and a low temperature seems
to be desirable.

Our choice for the risk estimate $\widetilde{r_t}$ is to use the
classical Stein unbiased estimate, which is sufficient to obtain optimal
oracle inequalities in expectation, 
\[ r_t=\|Y-\ft(Y)\|^2_2+ 2 \sigma^2 \tr(A_t)- n \sigma^2 \]
and add a penalty $\pen(t)$. We will consider simultaneously the case
of a penalty independent of $f_0$ and the one where the penalty may
depend on an upper bound of (kind of) sup norm.

More precisely, we allow the use, at least in the analysis, of an
upper bound $\widetilde{\|f_0-b\|}_{\infty}$ which can be thought as the
supremum of the sup norm of the coefficients of $f_0$ in any basis
appearing in $\mathcal{T}$. Indeed, we define
$\widetilde{\|f_0-b\|}_{\infty}$ as the smallest non-negative real
number $C$ such that for any $t \in \mathcal{T}$,
\begin{align*} 
\|A_t (f_0-b)\|_2^2 \leq C^2 \tr(A_t^2).
\end{align*}
By construction, $\widetilde{\|f_0-b\|}_{\infty}$ is smaller than the
sup norm of any coefficients of $f_0-b$ in any basis appearing in the
collection of estimators.
Note that $\widetilde{\|f_0-b\|}_{\infty}$ can also be
upper bounded by $\|f_0-b\|_1$, $\|f_0-b\|_2$ or
$\sqrt{n}\|f_0-b\|_{\infty}$ where the $\ell_1$ and sup norm can be
taken in any basis.

Our aim is to obtain sufficient conditions on the penalty $\pen(t)$
and the temperature $\beta$ so that an oracle inequality of type
\begin{align*}
  \|f_0 - f_{EWA}\|^2_2  \leq \inf_{\mu \in \mathcal{M}_+^1(\mathcal{T})}&  
(1 +\epsilon) \int \|f_0 - \ft\|^2_2 d\mu(t)\\ 
& \quad + (1 + \epsilon') \left( \int  \price(t) d\mu(t) + 2 \beta
   \KL(\mu,\pi) + \beta \ln\frac{1}{\eta} \right)
\end{align*}
holds either in probability or in expectation. Here, $\epsilon$ and
$\epsilon'$ are some small non-negative numbers
possibly equal to $0$ and $\price(t)$ a loss depending on the choice
of $\pen(t)$ and $\beta$. When $\mathcal{T}$ is countable, such an
oracle proves that the risk of our aggregate estimate is of the same
order as the one of the best estimate in the collection as it implies
\[
  \|f_0 - f_{EWA}\|^2_2  \leq \inf_{t \in \mathcal{T}}\left\lbrace  
(1 +\epsilon)  \|f_0 - \ft\|^2_2 
+ (1 + \epsilon') \left( \price(t) + \beta \ln\frac{1}{\pi(t)^2 \eta}
\right) \right\rbrace.
\]

Before stating our more general result, which is in
Section~\ref{sec:gener-oracle-ineq}, we provide a comparison with 
some similar results in the literature on the countable $\mathcal{T}$ setting.
	 
\section{Penalization strategies and preliminary results}
\label{sec:PenStrategies}

The most similar result in the literature is the one from
\citet{MR3192554} which holds under a Gaussian white noise assumption
and uses a penalty proportional to the known variance $\sigma^2$: 
\begin{prop}[\citet{MR3192554}]
  \label{prop:Dai}
  If $\pen(t)=2\sigma^2 \tr(A_t)$, and
  $\beta \geq 4\sigma^2 16$, then for all $\eta>0$, with
  probability at least $1-\eta$,
\begin{multline*} \|f_0-f_{EWA}\|^2 \leq \min_t \left\lbrace
\left(1+\frac{128 \sigma^2}{3\beta} \right) \|f_0-\ft\|^2
+8 \sigma^2 \tr(A_t)
\right. \\ \left.
+3\beta \ln \frac{1}{\pi_t} +3\beta \ln \frac{1}{\eta}
\right\rbrace.
\end{multline*} 
\end{prop}

Our result generalizes this result to the non necessarily independent
sub-Gaussian noise. We obtain
	\begin{prop}
	 \label{prop:simple_weak_gen}
	If $\beta\geq 20\sigma^2,$ there exists $\gamma \in [0,1/2)$, such that if \\$\pen(t) \geq \frac{4\sigma^2}{\beta-4\sigma^2} \tr(A_t^2) \sigma^2$, 
 for any $\eta>0$, with probability at least $1-\eta$,
\begin{align*} 
\|f_0-f_{EWA}\|^2 \leq \inf_{t} &\left\lbrace
\left(1+ \frac{4 \gamma}{1-2\gamma} \right) \|f_0-\ft\|^2 
\right. \\ 
&\left.+ \left( 1+ \frac{2\gamma}{1-2\gamma}\right) \left( \pen(t)+2\sigma^2 \tr(A_t) + 2\beta \ln
  \frac{1}{\pi_t} + \beta \ln \frac{1}{\eta}\right)
\right\rbrace.
\end{align*}
\end{prop}
The parameter $\gamma$  is 
explicit and satisfies $\epsilon =
O(\frac{\sigma^2}{\beta})$. We recover thus a similar weak oracle
inequality under a weaker assumption on the noise. It should be noted
that \cite{bellec} obtains a sharp oracle inequality for a slightly
different aggregation procedure but only under the very strong
assumption that $\tr(A_t) \leq \ln\frac{1}{\pi(t)}$.

Following \citet{MR3020421}, a lower bound on the
penalty, that involves the sup norm of $f_0$, can be given. In that case, the oracle
inequality is sharp as $\epsilon=\epsilon'=0$. Furthermore, the parameter $\gamma$ is not necessary and the minimum temperature is lower.
\begin{prop}
\label{prop:simple_sharp_gen}
If $\beta>4\sigma^2$, and
\[\pen(t)  \geq
    \frac{4\sigma^2}{\beta -4\sigma^2}\left( \sigma^2 \tr(A_t^2)+2
    \left[ \widetilde{\|f_0-b\|}_\infty^2 \tr(A_t^2) + \|b_t\|_2^2 \right] \right), \]  
then for any $\eta>0$, with probability at least $1-\eta$,
\begin{multline*} 
  \|f_0-f_{EWA}\|^2 \leq \inf_t \left\lbrace
    \|f_0-\ft\|^2 
    +2\sigma^2 
      \tr(A_t)
       \right. \\ \left.+\frac{8\sigma^2}{\beta-4\sigma^2}
      \left[ \widetilde{\|f_0-b\|}_\infty^2 \tr(A_t^2) + \|b_t\|_2^2 \right] 
  \right. \\ \left.  + \pen(t)+  2 \beta \ln \frac{1}{\pi_t} + \beta \ln \frac{1}{\eta}
\right\rbrace.
\end{multline*}
\end{prop}

The two results can be combined in a single one. Indeed, to obtain the
first oracle inequality, we rely in the proof on bounds of type
\begin{align*}
\|(A_t-A_u) f_0+b_t-b_u\|_2^2 \leq C_1 \|\ft-f_0\|_2^2+ C_2 \|\hat{f}_u -f_0\|_2^2, 
\end{align*}
with some constants $C_1$ and $C_2$ depending on $\gamma$ which  
allows to link $\|A_t f_0+b_t-A_u f_0+b_u \|_2^2$ to
$\|A_t Y+b_t-f_0 \|_2^2$ and $\|A_u Y+b_u-f_0 \|_2^2$. Whereas, for the second inequality
we rely on bounds of type
\begin{align*}
\|(A_t-A_u) f_0+b_t-b_u\|_2^2 &\leq 4 (\|A_t f_0\|_2^2+ \|A_u f_0\|_2^2+ \|b_t\|_2^2+\|b_u\|_2^2) \\
&\leq 4 \left[ \widetilde{\|f_0\|}_\infty^2 (\tr(A_t^2)+ \tr(A_u^2)) +\|b_t\|_2^2+\|b_u\|_2^2 \right]. 
\end{align*}	  
Combining these two upper bounds produce weak oracle inequalities for
a wider range of temperatures than
Proposition~\ref{prop:simple_weak_gen}, drawing a continuum between
Proposition~\ref{prop:simple_weak_gen} and
Proposition~\ref{prop:simple_sharp_gen}.  More precisely, one obtains
\begin{prop}
\label{prop:mixed-simple}
  For any $\delta \in [0,1]$, if $\beta \geq 4\sigma^2 V (1+4\delta)$
  and $\beta>4\sigma^2 V$, there exists $\gamma \geq 0$, such that if
  \[\pen(t)  \geq
    \frac{4\sigma^2}{\beta -4\sigma^2 V}\left( \sigma^2 \tr(A_t^2)+2(1-\delta)(1+2\gamma V)^2
    \left[ \widetilde{\|f_0\|}_\infty^2 \tr(A_t^2) + \|b_t\|_2^2 \right] \right), \]
  then for any $\eta>0$, with probability at least $1-\eta$,
  \begin{align*}
    \|f_0-f_{EWA}\|^2 \leq \inf_{t} &\left\lbrace \left(1+\epsilon
    \right) \|f_0-\ft\|^2 \right.\\
    &\left.+ (1 + \epsilon') \left( \price(t)
    +  2 \beta \ln \frac{1}{\pi_t} + \beta \ln \frac{1}{\eta}
\right) \right\rbrace.
  \end{align*}
with $\displaystyle \epsilon =   \frac{4V^2
          \gamma}{(2V-1)(1-2V\gamma)}$, $\displaystyle \epsilon'=
        \frac{2V\gamma}{1-2V\gamma}$ and
\begin{align*}
\price(t) & = \pen(t)+2\sigma^2
        \tr(A_t)+\frac{8\sigma^2(1-\delta)(1+2\gamma
          V)^2}{\beta-4\sigma^2 V}
        \left[ \widetilde{\|f_0\|}_\infty^2 \tr(A_t^2) + \|b_t\|_2^2 \right].
\end{align*}
\end{prop}
The convex combination parameter $\delta$ measures the account for
signal to noise ratio in the penalty.  We are now ready to state the
central result of this paper, which gives an explicit expression for
$\gamma$ and introduce an optimization parameter $\nu>0$.

\section{A general oracle inequality}
\label{sec:gener-oracle-ineq}

We consider now the general case for which $\mathcal{T}$ is not
necessarily countable. Recall that we have defined the exponentially weighted
measure $\rho_{EWA}$, a probability measure over $\mathcal{T}$, by
\begin{align*}
  d\rho_{EWA}(t) = \frac{\exp\left(-\frac{1}{\beta}\widetilde{r_t} \right)}
{\int \exp\left(-\frac{1}{\beta}\widetilde{r_{t'}} \right)
  d\pi(t')} d\pi(t)
\end{align*}
and the exponentially weighted aggregate $f_{EWA}$ by
$f_{EWA}=\int \ft \, d\rho_{EWA}(t)$. We will directly consider a
lower bound on the penalty of the same type than in
Proposition~\ref{prop:mixed-simple} and propositions similar to
Propositions~\ref{prop:simple_weak_gen}
and~\ref{prop:simple_sharp_gen} will be obtained as straightforward corollaries.

Our main contribution is the following two similar theorems:
\begin{thm}
  \label{thm:general_PACSimpleone}
  For any $\beta \geq 20 \sigma^2$, let
  \[ \gamma=\frac{\beta - 12 \sigma^2 -\sqrt{\beta -4
      \sigma^2}\sqrt{\beta -20 \sigma^2}}{16 \sigma^2} .  \]
  If for any $t \in \mathcal{T},$
  \[
    \pen(t)  \geq
    \frac{4\sigma^2}{\beta -4\sigma^2}  \sigma^2 \tr(A_t^2)  ,  
  \]
  then
  \begin{itemize}
  \item for any $\eta \in (0,1]$, with probability at least $1-\eta,$
    \begin{multline*}
      \|f_0-f_{EWA}\|_2^2 \leq \inf_{\mu \in
        \mathcal{M}_+^1(\mathcal{T})}
      \left(1+\frac{4 \gamma}{1-2\gamma} \right) \int \|f_0-\ft\|_2^2 d\mu(t)\\
      +\left(1+\frac{2 \gamma}{1-2\gamma}\right) \int  \pen(t) + 2\sigma^2 \tr(A_t)
      d\mu(t) +\beta \left(1+\frac{2 \gamma}{1-2\gamma}\right)
      \left(2\KL(\mu,\pi)+\ln\frac{1}{\eta} \right).
    \end{multline*}
  \item Furthermore
    \begin{multline*}
      \E \|f_0-f_{EWA}\|_2^2 \leq  \inf_{\mu \in
        \mathcal{M}_+^1(\mathcal{T})} \left(1+\frac{4 \gamma}{1-2\gamma} \right)
      \int \E
      \|f_0-\ft\|_2^2 d\mu(t)\\
      +\left(1+\frac{2 \gamma}{1-2\gamma}\right) \int  \pen(t) + 2\sigma^2 \tr(A_t)
      d\mu(t) +2\beta \left(1+\frac{2 \gamma}{1-2\gamma}\right)\KL(\mu,\pi).
    \end{multline*}
\end{itemize}
\end{thm}
and
\begin{thm}
  \label{thm:general_PACSimplezero}
  For any $\delta \in [0,1]$, if
  $\beta>4 \sigma^2$, 
  If for any $t \in \mathcal{T},$
  \[
    \pen(t)  \geq
    \frac{4\sigma^2}{\beta -4\sigma^2}\left( \sigma^2 \tr(A_t^2)+2
    \left[ \widetilde{\|f_0-b\|}_\infty^2 \tr(A_t^2) + \|b_t\|_2^2 \right] \right)  ,  
  \]
  then
  \begin{itemize}
  \item for any $\eta \in (0,1]$, with probability at least $1-\eta,$
    \begin{multline*}
      \|f_0-f_{EWA}\|_2^2 \leq \inf_{\mu \in
        \mathcal{M}_+^1(\mathcal{T})}
       \int \|f_0-\ft\|_2^2 d\mu(t)\\
      + \int \pen(t) + 2\sigma^2 \tr(A_t)+\frac{8\sigma^2}{\beta -4 \sigma^2} \left[\widetilde{\|f_0-b\|}_\infty^2 \tr(A_t^2)+\|b_t\|_2^2 \right] 
      d\mu(t)\\
      +\beta 
      \left(2\KL(\mu,\pi)+\ln\frac{1}{\eta} \right).
    \end{multline*}
  \item Furthermore
    \begin{multline*}
      \E \|f_0-f_{EWA}\|_2^2 \leq  \inf_{\mu \in
        \mathcal{M}_+^1(\mathcal{T})} \left(1+\frac{4 \gamma}{1-2\gamma} \right)
      \int \E
      \|f_0-\ft\|_2^2 d\mu(t)\\
      + \int \pen(t) + 2\sigma^2 \tr(A_t)+\frac{8\sigma^2}{\beta -4 \sigma^2} \left[\widetilde{\|f_0-b\|}_\infty^2 \tr(A_t^2)+\|b_t\|_2^2 \right] 
      d\mu(t) +2\beta  \KL(\mu,\pi).
    \end{multline*}
\end{itemize}
\end{thm}

When $\mathcal{T}$ is discrete, one can replace the minimization over
all the probability measure $\mathcal{M}^1_+(\mathcal{T})$ by the
minimization overall Dirac measure $\delta_{f_t}$ with $t\in\mathcal{T}$.  
Propositions~\ref{prop:simple_weak_gen} and~\ref{prop:simple_sharp_gen}
 are then straightforward corollaries. Note that the result in expectation
is obtained with the same penalty, which is known  not to be
necessary, at least in the Gaussian case, as shown by
\cite{MR3059085}.

If we assume the penalty is given
\begin{align*}
\pen(t) = \kappa \tr(A_t^2) \sigma^2,
\end{align*}
one can give rewrite the assumption in term of $\kappa$.
The weak oracle inequality holds for
any temperature greater than $20\sigma^2$ as soon as
$\kappa \geq \frac{4 \sigma^2}{\beta -4\sigma^2}$.
while an exact oracle inequality holds for any vector $f_0$ and any temperature $\beta$
greater than $4\sigma^2$ as soon as
\begin{align*}
\frac{\beta -4\sigma^2}{4\sigma^2 
}\kappa - 1 \geq
   \frac{\widetilde{\|f_0-b\|}_\infty^2 + \|b_t\|^2/\tr(A_t^2)}{\sigma^2}.
\end{align*}
For fixed $\kappa$ and $\beta$, this corresponds to a low peak signal
to noise ratio $\frac{\widetilde{\|f_0-b\|}_\infty^2}{\sigma^2}$ up to
the $\|b_t\|^2$ term which vanishes when $b_t=0$. Note that similar
results hold for a penalization scheme but with much larger constants
and some logarithmic factor in $n$.

Finally, the minimal temperature of $20\sigma^2$ can be replaced by
some smaller values if one further restrict the smoothed projections
used. As it appears in the proof, the temperature can be replaced by
$8\sigma^2$ or even $6\sigma^2$ when the smoothed
projections are respectively classical projections
and projections in the same basis. The
question of the minimality of such temperature is still open. Note
that in this proof, there is no loss due to the sub-Gaussianity
assumption, since the same upper bound on the exponential moment of
the deviation as in the Gaussian case are found, providing the same
penalty and bound on temperature.

The two results can be combined in a single one producing
 weak oracle inequalities for
a wider range of temperatures than
Theorem~\ref{thm:general_PACSimpleone}.
in Apprendix, we prove that a 
continuum between
those two cases exists: a weak oracle inequality, with smaller leading
constant than the one of Theorem~\ref{thm:general_PACSimpleone}, holds as
soon as there exists $\delta \in [0,1)$ such that
$\beta \geq 4 \sigma^2 (1+4\delta)$ and
\begin{align*}
\frac{\beta -4\sigma^2}{4\sigma^2 
}\kappa - 1 \geq
  (1-\delta)(1+2\gamma)^2 \frac{\widetilde{\|f_0-b\|}_\infty^2 + \|b_t\|^2/\tr(A_t^2)}{\sigma^2},
\end{align*}
where the signal to noise ratio guides the transition.  The
temperature required remains nevertheless always above $4 \sigma^2$.
The convex combination parameter $\delta$ measures the account for
signal to noise ratio in the penalty.

Note that in practice, the temperature can often be chosen smaller. It
is an open question whether the $4\sigma^2$ limit is an artifact of
the proof or a real lower bound. In the Gaussian case,
\cite{golubev14:_concen_inequal_expon_weigh_method} have been able to
show that this is mainly technical. Extending this result to our
setting is still an open challenge.

\appendix

\section{Proof of the oracle inequalities}
\label{sec:Proof_oracle}

The proof of this result is quite long and thus postponed in
Appendix~\ref{sec:Appendix_Proofs_subGaussian}. We provide first the
generic proof of the oracle inequalities, highlighting the role of
Gibbs measure and of some control in deviation. Then, we focus on the
aggregation of projection estimators in the Gaussian model. This
example already conveys all the ideas used in the complete proof of
the deviation lemma~: exponential moments inequalities for Gaussian
quadratic
form 
and the control of the bias $\|f_0-A_t f_0\|_2^2$ by
$\widetilde{\|f_0\|}_\infty^2$ on the one hand, to obtain an exact
oracle inequality, and by $\|f_0-A_t Y\|_2^2$ on the other hand,
giving a weak inequality.

The extension to the general case is obtained by showing that similar
exponential moments inequalities can be obtained for quadratic form of
sub-Gaussian random variables, working along the fact that the
systematic bias $\|f_0-A_t f_0\|_2^2$ is no longer always smaller than
$\|f_0-A_t Y\|_2^2$ and providing a fine tuning optimization allowing
the equality in the constraint on $\beta$ and an optimization on the
parameters $\epsilon$.

We provide in the next section the sketch of proof of
Theorem~\ref{thm:general_PAC}, an extended version of the Theorems as
well as its proof in the sub-Gaussian case and a simplified case
dealing with Gaussian noise and orthonormal projection  meant to be
compared with the one of \citet{MR3192554}.

\subsection{Extended result in the sub-Gaussian case} 
\label{sec:Appendix_Proofs_subGaussian}

We will consider affine estimators $\ft(Y)= A_t (Y-b) +
b + b_t$
corresponding to affine smoothed projection. We will assume that
\[ \ft(Y) = A_t (Y-b) + b + b_t =\sum_{i=1}^n \rho_{t,i} \langle Y-b, g_{t,i} \rangle
g_{t,i} + b + b_t \]
where $(g_{t,i})_{i=1}^n$ is an orthonormal basis,
$(\rho_{t,i})_{i=1}^n$ a sequence of non-negative real numbers and
$b_t \in \mathbb{R}^n$. By construction, $A_t$ is thus a symmetric
positive semi-definite real matrix. 
 We only assume here that the matrix collection
 $\{ A_t \}_{t \in \mathcal{T}}$ is such that
there exists a finite $V>0$ for which
$\sup_{t\in \mathcal{T}} \|A_t\|_2 \leq V$.
For sake of simplicity, we only use the
notation $\ft(Y)= A_t (Y-b) + b + b_t$ in the following. 

We obtain a theorem in which $V$ plays a role and in which a parameter
$\nu$ can be optimized.
\begin{thm}
  \label{thm:general_PAC}
  For any $\delta \in [0,1]$, if
  $\beta \geq 4 \sigma^2 V (1+4\delta ),$ and $\beta>4 \sigma^2 V$, let
  \[ \gamma=\frac{\beta -4 \sigma^2 V (1+2\delta)-\sqrt{\beta -4
      \sigma^2 V}\sqrt{\beta -4 \sigma^2 V (1+4\delta)}}{16 \sigma^2
    \delta V^2} \mathds{1}_{\delta>0}.  \]
  If for any $t \in \mathcal{T},$
  \[
    \pen(t)  \geq
    \frac{4\sigma^2}{\beta -4\sigma^2 V}\left( \sigma^2 \tr(A_t^2)+2(1-\delta)(1+2\gamma V)^2
    \left[ \widetilde{\|f_0-b\|}_\infty^2 \tr(A_t^2) + \|b_t\|_2^2 \right] \right)  ,  
  \]
  then
  \begin{itemize}
  \item for any $\eta \in (0,1]$, with probability at least $1-\eta,$
    \begin{multline*}
      \|f_0-f_{EWA}\|_2^2 \leq \inf_{\nu \in N} \inf_{\mu \in
        \mathcal{M}_+^1(\mathcal{T})}
      \left(1+\epsilon(\nu) \right) \int \|f_0-\ft\|_2^2 d\mu(t)\\
      +(1+\epsilon'(\nu)) \int \price(t)
      d\mu(t) +\beta (1+\epsilon'(\nu))
      \left(2\KL(\mu,\pi)+\ln\frac{1}{\eta} \right).
    \end{multline*}
  \item Furthermore
    \begin{multline*}
      \E \|f_0-f_{EWA}\|_2^2 \leq \inf_{\nu \in N} \inf_{\mu \in
        \mathcal{M}_+^1(\mathcal{T})} \left(1+\epsilon(\nu) \right)
      \int \E
      \|f_0-\ft\|_2^2 d\mu(t)\\
      +(1+\epsilon'(\nu)) \int \price(t)
      d\mu(t) +2\beta (1+\epsilon'(\nu))\KL(\mu,\pi),
    \end{multline*}
\end{itemize}
with $\displaystyle \epsilon(\nu)=\frac{1+\nu}{\nu}\frac{(1+\nu) \gamma}{1-(1+\nu)\gamma}$,
$\displaystyle \epsilon'(\nu)=\frac{(1+\nu)\gamma}{1-(1+\nu)\gamma}$,
\begin{align*}
\price(t)& = \pen(t) + 2\sigma^2 \tr(A_t)+\frac{8\sigma^2  (1-\delta)}{\beta -4 \sigma^2 V} (1+2\gamma
    V)^2 \left[\widetilde{\|f_0-b\|}_\infty^2 \tr(A_t^2)+\|b_t\|_2^2 \right]
\end{align*}
and $N=\{\nu>0 | (1+\nu)\gamma<1\}$.
\end{thm}

The parameter $\nu$ is a technical parameter that can be
optimized, provided $N$ is non empty. If $\delta > 0$, then 
for any $\beta \geq 4 \sigma^2 V(1+4\delta ),$ $0<2\gamma V\leq
1$. Thus $(0, 2V-1) \subseteq N$ as soon as $V>1/2$ with $2V-1 \in N$
if we assume that $\beta > 4 \sigma^2 V(1+4\delta )$. If we assume
$V\in(0,1/2)$, we have to impose $\beta > 4 \sigma^2 V + 2 \sigma^2
\delta (1+2V)^2$ in order to have a non empty $N$. Finally, if
$\delta=0$ then $\gamma=0$ and
$\epsilon'(\nu)=0, \epsilon(\nu)=0$, and no optimization is required.
Theorems~\ref{thm:general_PACSimpleone} and~\ref{thm:general_PACSimplezero} correspond to the case $V=1$ and
the choice $\nu= 2V-1=1$.

\subsection{General sketch of proof}
\label{subsec:General_Proof_oracle}

Theorem~\ref{thm:general_PAC} relies on the characterization of Gibbs
measure (Lemma~\ref{Gibbs}) and a control of deviation of the
empirical risk of any aggregate around its true risk.
			
$\rho$ is a Gibbs measure. Therefore it maximizes the entropy for a
given expected energy.  That is the subject of Lemma 1.1.3
in~\citet{MR2483528}: 
\begin{lem}
\label{Gibbs}
For any bounded measurable function
$h:\mathcal{T} \rightarrow \mathbb{R},$ and any probability
distribution $\rho \in \mathcal{M}_+^1\left(\mathcal{T}\right)$ such
that $\KL(\rho,\pi)<\infty,$
\[ \log \left( \int \exp(h) d\pi \right)= \int h d\rho-\KL(\rho,\pi)+\KL(\rho,\pi_{\exp(h)}),
 \]
 where by definition
 $\frac{d\pi_{\exp(h)}}{d\pi}=\frac{\exp[h(t)]}{\int \exp(h) d\pi }.$
 Consequently,
 \[ \log \left( \int \exp(h) d\pi \right)=\sup_{\rho \in
   \mathcal{M}_+^1\left(\mathcal{T}\right)} \int h
 d\rho-\KL(\rho,\pi). \]
\end{lem}

With $h(t)=-\frac{1}{\beta} [r_t+\pen(t)],$ this lemma states that for
any probability distribution
$\mu \in \mathcal{M}_+^1\left(\mathcal{T}\right)$ such that
$\KL(\mu,\pi)<\infty,$
\[  \int h d\rho -\KL(\rho,\pi) \geq \int h d\mu -\KL(\mu,\pi). \] 
Equivalently,
\begin{align*}
  &\int \|f_0-\ft\|_2^2 d\rho(t) +\int \left(r_t-\|f_0-\ft\|_2^2 +\pen(t) \right) d\rho(t) +\beta \KL(\rho,\pi)\\
  & \qquad \leq  
    \int \|f_0-\ft\|_2^2 d\mu(t) 
    +\int \left(r_t-\|f_0-\ft\|_2^2 +\pen(t) \right) d\mu(t) +\beta \KL(\mu,\pi)
  \\
  \Leftrightarrow
  &\int \|f_0-\ft\|_2^2 d\rho(t) -\int \|f_0-\ft\|_2^2 d\mu(t) \leq \int \left(\|f_0-\ft\|_2^2-r_t \right) d\rho(t)\\ 
  & \quad -\beta \KL(\rho,\pi) 
  -\int \left(\|f_0-\ft\|_2^2-r_t \right) d\mu(t)
    -\int \pen(t) d\rho(t)\\
   & \qquad +\int \pen(t)d\mu(t)
    +\beta \KL(\mu,\pi).
\end{align*}

The key is to upper bound the right-hand side with terms that may
depend on $\rho,$ but only through $\int \|f_0-\ft\|_2^2 d\rho(t)$ and
Kullback-Leibler distance. We will obtain two different controls in
the sub-Gaussian case and the Gaussian one that provide upper bounds
in probability (and in expectation) of type:
\begin{multline*}
  \int\left(\|f_0-\ft\|^2_2-r_t\right)d\rho(t)-\int\left(\|f_0-\hat{f}_u\|^2_2-r_u
  \right) d\mu(u)
  \\
  \leq C_1 \int \|f_0-\ft\|_2^2 d\rho(t) + C_2 \int
  \|f_0-\hat{f}_u\|_2^2 d\mu(u)
  \\
  + \int \left(C_3\tr(A_t^2)+ C_4 \|b_t\|_2^2\right) d\rho(t) \\
  +  C_5 \int
  \tr(A_u) d\mu(u) 
  + \int \left(C_6\tr(A_u^2)+ C_7 \|b_u\|_2^2\right) d\mu(u) 
  \\
  +\beta \left( \KL(\rho,\pi)+ \KL(\mu,\pi)+ \ln\frac{1}{\eta} \right)
\end{multline*}
where $C_1$ to $C_7$ are known functions.  Combining with the previous
inequality and taking $\pen(t) \geq C_3 \tr(A_t^2)+ C_4 \|b_t\|_2^2$ gives
\begin{multline*}
  (1-C_1)\int \|f_0-\ft\|_2^2 d\rho(t) -(1+C_2)\int \|f_0-\ft\|_2^2 d\mu(t) \\
  \leq C_5 \int \tr(A_u) d\mu(u) + \int \left(C_6\tr(A_u^2)+ C_7 \|b_u\|_2^2\right) d\mu(u) +\int
  \pen(u)d\mu(u)
  \\
  +\beta \left( 2 \KL(\mu,\pi)+ \ln\frac{1}{\eta} \right).
\end{multline*}
The additional condition $C_1<1$ allows to conclude.  It is now clear
that the whole work lies in the proof of the lemma.

\subsection{Proof of Theorem~\ref{thm:general_PAC}}
\label{subsec:Appendix_Proofs_oracle}

		The proof follows from the scheme described in section~\ref{subsec:General_Proof_oracle}. The main point is to control \[\int \left(\|f_0-\ft\|_2^2-r_t \right) d\rho(t)-\int \left(\|f_0-\ft\|_2^2-r_t \right) d\mu(t). \] We recall that $A_t$ is a symmetric positive semi-definite matrix, there exists $V>0$ such that $\sup_{t \in \mathcal{T}} \|A_t\|_2 \leq V$ and $W$ is a centered sub-Gaussian noise. 
For any $t, u \in \mathcal{T},$ we denote $\Delta_{t,u}=\|f_0-\ft\|^2_2-r_t-\|f_0-\hat{f}_u\|^2_2+r_u.$

The exponential moment of $\Delta_{t,u}$ is easily controlled by a
term involving $\|(A_t-A_u)(f_0-b)+b_t-b_u\|_2^2$ (see
Equation~\eqref{eq:Delta_t,u}). 
In the projection case, we used a bound
\begin{align*}
\|(A_t-A_u)(f_0-b)\|_2^2  
&\leq 2 \left( \|A_t(f_0-b)
  - (f_0-b) \|^2 + \|A_u(f_0-b)
  -(f_0-b)\|^2 \right) \\
&\leq 2 \left(\|A_tY - (f_0-b)\|^2 +\|A_uY-(f_0-b) \|^2 \right) 
\end{align*}
whose last inequalities  is not applicable here.
To overcome this difficulty, a term $\|(A_t-A_u)Y\|_2^2$ is introduced and for an arbitrary $\gamma \geq 0$, we try to control $\Delta_{t,u}-\gamma \|(A_t-A_u)Y\|_2^2$.

A simple calculation yields
\begin{multline*}
\Delta_{t,u}-\gamma \|\ft-\hat{f}_u\|_2^2=
W^\top (2I-\gamma (A_t-A_u)^\top)(A_t-A_u)W \\
+ 2 W^\top (I-\gamma (A_t-A_u)^\top)\left[(A_t-A_u)(f_0-b)+b_t-b_u \right]\\
- 2\sigma^2 \tr(A_t-A_u)-\gamma \|(A_t-A_u)(f_0-b)+b_t-b_u\|_2^2 .
\end{multline*}
Using $W^\top (2I-\gamma (A_t-A_u)^\top)(A_t-A_u)W \leq 2 W^\top (A_t-A_u)W$ and since $(A_t)_{t \in \mathcal{T}}$ are positive semi-definite matrices, $2 W^\top(A_t-A_u)W \leq 2 W^\top A_t W$. Thus, for any $\beta>0,$ any $\gamma \geq 0$, 
\begin{multline*}
\E\exp \left(\frac{\Delta_{t,u}}{\beta}-\frac{\gamma}{\beta} \|\ft-\hat{f}_u\|_2^2\right)\\
\leq
\E\left[\exp \frac{2}{\beta}\left( W^\top A_t W +W^\top (I-\gamma (A_t-A_u))\left[(A_t-A_u)(f_0-b)+b_t-b_u \right] \right)\right]\\
\times \exp \frac{-1}{\beta}\left( 2 \sigma^2 \tr(A_t-A_u)+ \gamma \|(A_t-A_u)(f_0-b)+b_t-b_u\|_2^2 \right)  .
\end{multline*}

 	 The first step is to bring us back to the Gaussian case, using $W$'s sub-Gaussianity and an idea of \citet{MR2994877}. Let $Z$ be a standard Gaussian random variable, independent of $W$. Then,
\begin{align*}
&\E \exp \left( \frac{2}{\sqrt{\beta}} W^\top \sqrt{A_t} Z+\frac{2}{\beta} W^\top (I-\gamma (A_t-A_u))\left[(A_t-A_u)(f_0-b)+b_t-b_u \right] \right) \\
& \quad= \E \left[ \E \left[ \exp \left( \frac{2}{\sqrt{\beta}} W^\top \sqrt{A_t} Z+\frac{2}{\beta} W^\top (I-\gamma (A_t-A_u))\left[(A_t-A_u)(f_0-b)+b_t-b_u \right] \right) \big|W \right] \right]\\
& \quad=\E \left[ \E \left[ \exp \left( \frac{2}{\sqrt{\beta}} W^\top \sqrt{A_t} Z  \right) \big|W \right] 
\right. \\ 
& \left. \qquad \qquad \qquad \qquad \qquad   \times
\exp \left(\frac{2}{\beta} W^\top (I-\gamma (A_t-A_u))\left[(A_t-A_u)(f_0-b)+b_t-b_u \right] \right)  \right] \\
& \quad=\E \exp \frac{2}{\beta} \left(W^\top A_t W+ W^\top (I-\gamma (A_t-A_u))\left[(A_t-A_u)(f_0-b)+b_t-b_u \right] \right).
\end{align*}
On the other hand,
\begin{multline*}
\E\left[\exp \frac{2}{\beta}\left( W^\top A_t W +W^\top (I-\gamma (A_t-A_u))\left[(A_t-A_u)(f_0-b)+b_t-b_u \right] \right)\right]\\
= \E \left[ \E \left[ \exp \left( \frac{2}{\sqrt{\beta}} W^\top \sqrt{A_t} Z+\frac{2}{\beta} W^\top (I-\gamma (A_t-A_u))\left[(A_t-A_u)(f_0-b)+b_t-b_u \right] \right) \big|Z \right] \right].
\end{multline*}
Since $W$ is sub-Gaussian with parameter $\sigma$,
\begin{multline*}
\E\left[\exp \frac{2}{\beta}\left( W^\top A_t W +W^\top (I-\gamma (A_t-A_u))\left[(A_t-A_u)(f_0-b)+b_t-b_u \right] \right)\right]\\
\leq \E \exp \left( \frac{\sigma^2}{2} \left \| \frac{2}{\sqrt{\beta}} \left(\sqrt{A_t}Z +\frac{1}{\sqrt{\beta}} (I-\gamma (A_t-A_u))\left[(A_t-A_u)(f_0-b)+b_t-b_u \right] \right)\right \|_2^2 \right)
\end{multline*}
Hence,
\begin{multline*}
\E\exp \left(\frac{\Delta_{t,u}}{\beta}-\frac{\gamma}{\beta} \|\ft-\hat{f}_u\|_2^2\right)\\
\leq
\E\left[\exp \frac{2\sigma^2}{\beta}\left( Z^\top A_t Z+\frac{2}{\sqrt{\beta}} Z^\top \sqrt{A_t} (I-\gamma (A_t-A_u))\left[(A_t-A_u)(f_0-b)+b_t-b_u \right] \right)\right]\\
\times \exp \left( \frac{2\sigma^2}{\beta^2} \left \|(I-\gamma (A_t-A_u))\left[(A_t-A_u)f_0+b_t-b_u \right] \right \|_2^2-\frac{2 \sigma^2}{\beta} \tr(A_t-A_u)\right)\\
\times \exp \left(-\frac{\gamma}{\beta} \|(A_t-A_u)(f_0-b)+b_t-b_u\|_2^2 \right) .
\end{multline*}
 The expectation is the one of the exponential of some quadratic form
 and we will use the ideas of~\citet{MR2994877}. Since $A_t$ is positive semi-definite, there exist an orthogonal matrix $U$ and a diagonal matrix $D$ such that $A_t=U^\top D U.$	Note that $UZ$ is a standard Gaussian variable.
 This diagonalization step and the non-negativity of the eigenvalues allow to apply Lemma 2.4 of~\citet{MR2994877}. Then, for any $\beta > 4\sigma^2 V$, any $\gamma \geq 0$,
 \begin{multline*}
\E\exp \left(\frac{\Delta_{t,u}}{\beta}-\frac{\gamma}{\beta} \|\ft-\hat{f}_u\|_2^2\right)
\leq
\exp \left(\frac{2\sigma^2}{\beta} \tr(A_t) + \frac{4\sigma^4}{\beta (\beta-4\sigma^2 V)}  \tr(A_t^2)\right)\\
\times \exp \left( \frac{8\sigma^4}{\beta^2 (\beta-4\sigma^2 V)}  \left \| \sqrt{A_t} (I-\gamma (A_t-A_u))\left[(A_t-A_u)(f_0-b)+b_t-b_u \right] \right \|_2^2 \right) \\
\times \exp \left( \frac{2\sigma^2}{\beta^2 } \left \|(I-\gamma (A_t-A_u))\left[(A_t-A_u)(f_0-b)+b_t-b_u \right] \right \|_2^2-\frac{2 \sigma^2}{\beta} \tr(A_t-A_u)\right)\\
\times \exp\left(-\frac{\gamma}{\beta} \|(A_t-A_u)(f_0-b)+b_t-b_u\|_2^2 \right) .
\end{multline*}
Consequently,
 \begin{align*}
&\E\exp \left(\frac{\Delta_{t,u}}{\beta}+\frac{\gamma}{\beta} \left(\|(A_t-A_u)(f_0-b)+b_t-b_u\|_2^2-\|\ft-\hat{f}_u\|_2^2 \right)\right)\\
&\leq
\exp \frac{2\sigma^2}{\beta}\left( \tr(A_u)+\frac{2\sigma^2}{\beta-4\sigma^2 V} \tr(A_t^2) \right)\\
&\quad \times \exp \left( \frac{2\sigma^2}{\beta^2} \left(\frac{4\sigma^2 V}{\beta-4 \sigma^2 V}(1+2\gamma V)^2+(1+2\gamma V)^2 \right) \|(A_t-A_u)(f_0-b)+b_t-b_u\|_2^2\right)\\
&\leq \exp \frac{2\sigma^2}{\beta}\left( \tr(A_u)+\frac{2\sigma^2}{\beta-4\sigma^2 V} \tr(A_t^2) +\frac{(1+2\gamma V)^2}{\beta -4 \sigma^2 V}  \|(A_t-A_u)(f_0-b)+b_t-b_u\|_2^2\right) .
\end{align*}

If an exact oracle inequality is wished, $\|(A_t-A_u)(f_0-b)+b_t-b_u\|_2^2$ should be upper bounded by some constant and $\gamma$ should be set to zero. Else, $\gamma$ is used to \emph{replace} the terms in $\|(A_t-A_u)(f_0-b)+b_t-b_u\|_2^2$ by $\|(A_t-A_u)(Y-b)+b_t-b_u\|_2^2$.
Thus, the terms depending on $f_0$ will be upper bounded in two ways: \begin{itemize}
\item on the one hand, using $\widetilde{\|f_0-b\|}_\infty^2$
\begin{multline*}
\|(A_t-A_u)(f_0-b)+b_t-b_u\|_2^2 
\leq 4 \left(\|A_t (f_0-b)\|_2^2+\|A_u (f_0-b)\|_2^2 +\|b_t\|_2^2+\|b_u\|_2^2\right) \\
\leq 4 \left(\left( \tr(A_t^2)+ \tr(A_u^2) \right) \widetilde{\|f_0-b\|}_\infty^2+\|b_t\|_2^2+\|b_u\|_2^2\right)
\end{multline*}

For any $\delta \in [0,1]$,
 \begin{align*}
&\E\exp \left(\frac{\Delta_{t,u}}{\beta}+\frac{\gamma}{\beta} \left(\|(A_t-A_u)(f_0-b)+b_t-b_u\|_2^2-\|\ft-\hat{f}_u\|_2^2 \right)\right)\\
&\leq \exp \frac{2\sigma^2}{\beta}\left( \tr(A_u)+\frac{2\sigma^2}{\beta-4\sigma^2 V} \tr(A_t^2) \right) \\
&\qquad  \times \exp \left( \frac{2\sigma^2}{\beta} \frac{ (1+2\gamma V)^2 (1-\delta)}{\beta -4 \sigma^2 V}  \|(A_t-A_u)(f_0-b)+b_t-b_u\|_2^2\right)\\
&\qquad \qquad \times \exp \left(\frac{2\sigma^2 (1+2\gamma V)^2 \delta}{\beta(\beta -4 \sigma^2 V)}  \|(A_t-A_u)(f_0-b)+b_t-b_u\|_2^2 \right)\\
&\leq \exp \frac{2\sigma^2}{\beta}\left( \tr(A_u)+\frac{2\sigma^2}{\beta-4\sigma^2 V} \tr(A_t^2)+\frac{ (1+2\gamma V)^2 \delta}{\beta -4 \sigma^2 V}  \|(A_t-A_u)(f_0-b)+b_t-b_u\|_2^2  \right)\\
&\times \exp\left(\frac{8\sigma^2 (1+2\gamma V)^2 (1-\delta)}{\beta(\beta-4 \sigma^2 V)}  \left[ \left(\tr(A_t^2)+\tr(A_u^2)\right) \widetilde{\|f_0-b\|}_\infty^2 +\|b_t\|_2^2+\|b_u\|_2^2\right] \right).
\end{align*}

\item on the other hand, introducing $\|\ft-f_0\|_2^2$ to obtain a weak oracle inequality: conditions should be found on $\gamma$ such that 
\begin{multline*}
\frac{2\sigma^2 (1+2\gamma V)^2 \delta}{\beta-4 \sigma^2 V}  \|(A_t-A_u)(f_0-b)+b_t-b_u\|_2^2 \\
-\gamma \left(\|(A_t-A_u)(f_0-b)+b_t-b_u\|_2^2-\|\ft-\hat{f}_u\|_2^2\right)\\
 \leq C_1 \|\ft-f_0\|_2^2+ C_2 \|\hat{f}_u-f_0\|_2^2
\end{multline*}
for some non-negative constants $C_1$ and $C_2$ and with $\delta>0$.
Since for any $\nu>0$, $\|\ft-\hat{f}_u\|_2^2 \leq (1+\nu) \|\ft-f_0\|_2^2+ \left(1+\frac{1}{\nu}\right) \|\hat{f}_u-f_0\|_2^2$, it suffices that \[
\frac{2\sigma^2 (1+2\gamma V)^2 \delta}{\beta -4 \sigma^2 V}  \|(A_t-A_u)(f_0-b)+b_t-b_u\|_2^2 -\gamma \|(A_t-A_u)(f_0-b)+b_t-b_u\|_2^2 \leq 0.
\]
This condition may be fulfilled if $\beta \geq 4\sigma^2 V(1+4\delta)$. The smallest $\gamma\geq 0$ among all the possible ones is chosen~:
\[ \gamma= \frac{1}{16\sigma^2 \delta V^2} \left(\beta -4\sigma^2 V (1+2\delta)-\sqrt{\beta -4\sigma^2 V}\sqrt{\beta -4\sigma^2 V(1+4\delta)} \right) \mathds{1}_{\delta>0}. \]
\end{itemize}

This leads to the following inequality~: for any $\delta \in [0,1]$,
for any $\beta >4 \sigma^2 V$ and $\beta \geq 4 \sigma^2 V (1+4\delta)$, with $\gamma$ previously defined, for any $\nu>0$,
 \begin{multline*}
\E\exp \left(\frac{\Delta_{t,u}}{\beta}
 -\frac{\gamma}{\beta}\left((1+\nu)\|\ft-f_0\|_2^2 +\left(1+\frac{1}{\nu}\right)\|\hat{f}_u-f_0\|_2^2 \right)\right)\\
\leq \exp\left(\frac{8\sigma^2 (1+2\gamma V)^2 (1-\delta)}{\beta(\beta-4 \sigma^2 V)}  \left[ \left(\tr(A_t^2)+\tr(A_u^2)\right) \widetilde{\|f_0-b\|}_\infty^2 +\|b_t\|_2^2+\|b_u\|_2^2\right] \right)\\
\times \exp \frac{2\sigma^2}{\beta}\left( \tr(A_u)+\frac{2\sigma^2}{\beta-4\sigma^2 V} \tr(A_t^2) \right).
\end{multline*}

Along the same lines as~\citet{MR2786484}, we first integrate according to the prior $\pi$, use Fubini's theorem, introduce the probability measures $\rho$ and $\mu$ and apply Jensen's inequality to obtain that for any $\eta \in (0,1]$,
 \begin{multline}
 \label{eq:Delta_t,u_General}
\E \exp \frac{1}{\beta} \Big[ \int \int \Delta_{t,u} d\rho(t) d\mu(u) 
-(1+\nu)\gamma \int \|\ft-f_0\|_2^2 d\rho(t)\\
-\frac{4\sigma^2}{\beta -4\sigma^2 V} \int \left( \sigma^2 \tr(A_t^2)+2(1-\delta)(1+2\gamma V)^2
    \left[ \widetilde{\|f_0-b\|}_\infty^2 \tr(A_t^2) + \|b_t\|_2^2 \right] \right)d\rho(t)\\
-2\sigma^2 \left( \int \tr(A_u) d\mu(u) +\frac{4(1-\delta)(1+2 \gamma V)^2}{\beta -4 \sigma^2 V} \int \left[ \widetilde{\|f_0-b\|}_\infty^2 \tr(A_u^2) + \|b_u\|_2^2 \right] d\mu(u) \right)\\
-\left(1+\frac{1}{\nu}\right)\gamma \int \|\hat{f}_u-f_0\|_2^2 d\mu(u)
-\beta \left( \KL(\rho,\pi)+\KL(\mu,\pi)+\ln\frac{1}{\eta}\right)
\Big] \leq \eta.
 \end{multline}
 Finally, using $\exp(x)\geq \mathds{1}_{\mathbb{R}_+}(x),$ for any $\delta \in [0,1]$, any $\beta>4\sigma^2 V$ and $\beta \geq 4\sigma^2 V (1+4\delta)$, with $\gamma$ previously defined, 
 for any $\eta \in (0,1]$, for any $\nu >0$,
\begin{align*}
&\mathbb{P}\Big[ 
\int \int \Delta_{t,u} d\rho(t) d\mu(u) 
\leq
(1+\nu)\gamma \int \|\ft-f_0\|_2^2 d\rho(t)\\
&+\frac{4\sigma^2}{\beta -4\sigma^2 V} \int \left( \sigma^2 \tr(A_t^2)+2(1-\delta)(1+2\gamma V)^2
\left[ \widetilde{\|f_0-b\|}_\infty^2 \tr(A_t^2) + \|b_t\|_2^2 \right] \right)d\rho(t)\\
&+2\sigma^2 \left( \int \tr(A_u) d\mu(u) +\frac{4(1-\delta)(1+2 \gamma V)^2}{\beta -4 \sigma^2 V} \int \left[ \widetilde{\|f_0-b\|}_\infty^2 \tr(A_u^2) + \|b_u\|_2^2 \right] d\mu(u) \right)\\
&+\left(1+\frac{1}{\nu}\right)\gamma \int \|\hat{f}_u-f_0\|_2^2 d\mu(u)
+\beta \left( \KL(\rho,\pi)+\KL(\mu,\pi)+\ln\frac{1}{\eta}\right)
\Big] \geq 1-\eta.
\end{align*}

Combining Equation~\eqref{eq:Delta_t,u_General} with $\eta=1$ and the
inequality $t \leq \exp(t)-1$ leads to
\begin{multline*}
\E \left[ \int \int \Delta_{t,u} d\rho(t) d\mu(u) \right]
\leq \E\Big[
(1+\nu)\gamma \int \|\ft-f_0\|_2^2 d\rho(t)\\
+\frac{8\sigma^2}{\beta -4\sigma^2 V} (1-\delta)(1+2\gamma V)^2 \int
\left[ \widetilde{\|f_0-b\|}_\infty^2 \tr(A_t^2) + \|b_t\|_2^2 \right]  d\rho(t)\\
+\frac{4\sigma^4}{\beta -4\sigma^2 V} \int  \tr(A_t^2) d\rho(t)
+2\sigma^2 \int \tr(A_u) d\mu(u) \\
+ \frac{8\sigma^2}{\beta -4 \sigma^2 V} (1-\delta)(1+2\gamma V)^2 \int \left[\widetilde{\|f_0-b\|}_\infty^2 \tr(A_u^2) +\|b_u\|_2^2 \right] d\mu(u) \\
+\left(1+\frac{1}{\nu}\right)\gamma \int \|\hat{f}_u-f_0\|_2^2 d\mu(u)
+\beta \left( \KL(\rho,\pi)+\KL(\mu,\pi) \right)\Big].
\end{multline*}

We deduce thus that with probability at least $1-\eta$,
\begin{multline*}
\int \|f_0-\ft\|_2^2 d\rho(t) -\int \|f_0-\ft\|_2^2 d\mu(t) 
\leq
(1+\nu)\gamma \int \|\ft-f_0\|_2^2 d\rho(t)\\
+\frac{4\sigma^2}{\beta -4\sigma^2 V} \int \left( \sigma^2 \tr(A_t^2)+2(1-\delta)(1+2\gamma V)^2
\left[ \widetilde{\|f_0-b\|}_\infty^2 \tr(A_t^2) + \|b_t\|_2^2 \right] \right)  d\rho(t) \\
+2\sigma^2 \left( \int \tr(A_t) d\mu(t) +\frac{4(1-\delta)(1+2\gamma V)^2}{\beta -4 \sigma^2 V}  \int \left[\widetilde{\|f_0-b\|}_\infty^2 \tr(A_t^2)+\|b_t\|_2^2 \right] d\mu(t) \right)\\
-\int \pen(t) d\rho(t)+\int \pen(t)d\mu(t)
+\left(1+\frac{1}{\nu}\right)\gamma \int \|\ft-f_0\|_2^2 d\mu(t)\\
+\beta \left(2 \KL(\mu,\pi)+\ln\frac{1}{\eta}\right).
 \end{multline*}
Taking \[
    \pen(t)  \geq
    \frac{4\sigma^2}{\beta -4\sigma^2 V}\left( \sigma^2 \tr(A_t^2)+2(1-\delta)(1+2\gamma V)^2
    \left[ \widetilde{\|f_0-b\|}_\infty^2 \tr(A_t^2) + \|b_t\|_2^2 \right] \right)  ,  
  \]and $\nu \in N=\{\nu>0|(1+\nu)\gamma<1\}$, such that the inequality stays informative,
\begin{multline*}
\left(1-(1+\nu)\gamma\right) \int \|f_0-\ft\|_2^2 d\rho(t)
\leq \left(1+\left(1+\frac{1}{\nu}\right)\gamma \right) \int \|f_0-\ft\|_2^2 d\mu(t) \\
+2\sigma^2 \left( \int \tr(A_t) d\mu(t) +\frac{4(1-\delta)(1+2\gamma V)^2}{\beta -4 \sigma^2 V}  \int \left[\widetilde{\|f_0-b\|}_\infty^2 \tr(A_t^2)+\|b_t\|_2^2 \right] d\mu(t) \right)\\
+\int \pen(t)d\mu(t)+\beta \left(2 \KL(\mu,\pi)+\ln\frac{1}{\eta}\right).
 \end{multline*}
Finally, since $\|f_0-f_{EWA}\|_2^2 \leq \int \|f_0-\ft\|_2^2 d\rho(t)$, 
\begin{multline*}
\|f_0-f_{EWA}\|_2^2 
\leq \left(1+\frac{(1+\nu)^2\gamma}{\nu(1-(1+\nu)\gamma)} \right) \int \|f_0-\ft\|_2^2 d\mu(t) \\
+\frac{2\sigma^2}{1-(1+\nu)\gamma} \left( \int \tr(A_t) d\mu(t) 
\right.\\
\left. 
+\frac{4(1-\delta)(1+2\gamma V)^2}{\beta -4 \sigma^2 V}  \int \left[\widetilde{\|f_0-b\|}_\infty^2 \tr(A_t^2)+\|b_t\|_2^2 \right] d\mu(t) \right)\\
+\frac{1}{1-(1+\nu)\gamma} \left(\int \pen(t)d\mu(t)+\beta \left(2 \KL(\mu,\pi)+\ln\frac{1}{\eta}\right) \right).
 \end{multline*}
The result in expectation is obtained in the same fashion.

\subsection{Gaussian noise case and projection estimates}
\label{subsec:Gaussian_proof}

In this subsection,  we assume that
$A_t$ are the matrices of orthogonal projections, $b_t=0$, and the
noise $W$ is a centered Gaussian random variable with variance
$\sigma^2 I$.  The constants in the previous theorem can be enhanced:
\begin{thm}
  \label{thm:Gaussian}
  Let $\pi$ be an arbitrary prior measure over $\mathcal{T}$. For any
  $\delta \in [0,1]$, any $\beta>4\sigma^2(\delta+1)$, the aggregate
  estimator $f_{EWA}$ defined with
  \[ \pen(t) \geq \frac{2\sigma^4}{\beta -4\sigma^2}
    \left(1+2(1-\delta)
      \frac{\widetilde{\|f_0-b\|}_\infty^2}{\sigma^2} \right)
    \tr(A_t) \] satisfies the oracle inequalities of
  Theorem~\ref{thm:general_PAC} with
  $\displaystyle \epsilon = 2 \epsilon' = \frac{8\sigma^2
    \delta}{\beta -4\sigma^2 (\delta+1)}$ and
  \[
    \price(t) = \pen(t)+2\left(1+\frac{2(1-\delta)\sigma^2}{\beta -4
        \sigma^2} \frac{\widetilde{\|f_0-b\|}_\infty^2}{\sigma^2}
    \right) \tr(A_t)\sigma^2
  \]
\end{thm}

\noindent Note that the result may be further simplified using $\price(t) \leq 2
\left( \pen(t)+ \sigma^2 \tr(A_t)\right)$.

Again, the key is a control of the
deviation of the empirical risk of any aggregate around its true
risk.  We focus now on the proof of such a control obtaied by mixing
control of exponential moments of a quadratic form of a Gaussian
random variable with basic inequalities like Jensen, Fubini, and the
important link between $\|f_0-A_t f_0\|_2^2$ and $\|f_0-A_t
Y\|_2^2$. 
For the sake of clarity, for any $t, u \in \mathcal{T},$ let \[\Delta_{t,u}=\|f_0-\ft\|^2_2-r_t-\|f_0-\hat{f}_u\|^2_2+r_u.\]
A simple calculation yields
\[
\Delta_{t,u}=2\left(W^\top(A_t-A_u)W +W^\top(A_t-A_u)(f_0-b)-\sigma^2 \tr(A_t-A_u) \right).
\]
Since $(A_t)_{t \in \mathcal{T}}$ are positive semi-definite matrices, $W^\top(A_t-A_u)W \leq W^\top A_t W,$ and there exist an orthogonal matrix $U$ and a diagonal matrix $D$ such that $A_t=U^\top D U.$

For any $\beta>0,$ 
\begin{multline*}
\E\left[\exp \frac{\Delta_{t,u}}{\beta}\right]
\leq
\E\left[\exp \frac{2}{\beta}\left( (UW)^\top D(UW)
 +(UW)^\top U(A_t-A_u)(f_0-b)
\right. \right. \\ \left. \left.
-\sigma^2 \tr(A_t-A_u)\right) \right] .
\end{multline*}

Following lemma 2.4 of~\citet{MR2994877}, if $\beta>4 \sigma^2,$ 
\begin{align}
 \label{eq:Delta_t,u}
\E\left[\exp \frac{\Delta_{t,u}}{\beta}\right]
\leq
\exp \frac{2\sigma^2}{\beta}\left(\tr(A_u)+\frac{2\sigma^2 \tr(A_t)+ \|(A_t-A_u)(f_0-b)\|_2^2 }{\beta-4\sigma^2}\right).
\end{align}
Note that 
\begin{align*}
\|(A_t-A_u)(f_0-b)\|_2^2 &\leq 2 \left(\|(f_0-b) -A_t (f_0-b) \|_2^2+
  \|(f_0-b) - A_u (f_0-b) \|_2^2 \right)\\
&
\leq 2 \left(\|(f_0 - b) - A_t (Y-b)\|_2^2+\|(f_0-b) - A_u (Y-b)\|_2^2
  \right)\\
& \leq 2 \left(\|f_0 - \ft\|^2 + \|f_0 - \hat{f}_u\|^2 \right)
\end{align*}
and 
\[ 
\|(A_t-A_u)(f_0-b)\|_2^2 \leq 2 \left(\|A_t (f_0-b)\|_2^2+\|A_u (f_0-b)\|_2^2 \right)
\leq 2 \widetilde{\|f_0-b\|}_\infty^2 \left( \tr(A_t)+\tr(A_u)\right).
\]
Thus, for any $\beta>4\sigma^2,$ for any $\delta \in [0,1],$
\begin{multline*}
\E \exp \left[\frac{\Delta_{t,u}}{\beta}
- \frac{2\sigma^2}{\beta}\left(\tr(A_u)+\frac{2\sigma^2 \tr(A_t)}{\beta -4\sigma^2} \right) \right. \\
\left.
-\frac{4\sigma^2 \delta}{\beta (\beta -4\sigma^2)}
 \left(\|f_0-\ft\|_2^2+\|f_0-\hat{f}_u\|_2^2 \right)
\right. \\
\left.
-\frac{4\sigma^2}{\beta  (\beta -4\sigma^2)}
 (1-\delta) \widetilde{\|f_0-b\|}_\infty^2 \left( \tr(A_t)+\tr(A_u)\right)
 \right]\leq 1.
\end{multline*}

The proof can be concluded, along the same lines as~\citet{MR2786484}, by first integrating according to the prior $\pi$, using Fubini's theorem, introducing the probability measures $\rho$ and $\mu$ and applying Jensen's inequality.

\bibliographystyle{plainnat}

\end{document}